\documentclass[reqno]{amsart}

\usepackage{epsf,amssymb,graphics,graphicx,wrapfig}

\def\eps{\varepsilon}

\def\be{\begin{equation}}
\def\ee{\end{equation}}
\def\ba{\begin{align}}
\def\bm{\begin{multline}}
\def\bfig{\begin{figure}[htb]}
\def\efig{\end{figure}}

\numberwithin{equation}{section}
\newtheorem{theorem}{Theorem}[section]
\newtheorem{proposition}[theorem]{Proposition}
\newtheorem{lemma}[theorem]{Lemma}

\newcommand{\nn}{\nonumber}

\DeclareMathSymbol{\leqslant}{\mathalpha}{AMSa}{"36}
\DeclareMathSymbol{\geqslant}{\mathalpha}{AMSa}{"3E}
\DeclareMathSymbol{\doteqdot}{\mathalpha}{AMSa}{"2B}
\DeclareMathSymbol{\circlearrowright}{\mathalpha}{AMSa}{"08}
\DeclareMathSymbol{\subsetneq}{\mathalpha}{AMSb}{"28}
\DeclareMathSymbol{\supsetneq}{\mathalpha}{AMSb}{"29}
\renewcommand{\leq}{\;\leqslant\;}
\renewcommand{\geq}{\;\geqslant\;}

\newcommand{\dd}{{\rm d}}
\newcommand{\e}[1]{\,{\rm e}^{#1}\,}
\newcommand{\ii}{{\rm i}}
\newcommand{\sumtwo}[2]{\sum_{\substack{#1 \\ #2}}}

\newcommand{\suptwo}[2]{\sup_{\substack{#1 \\ #2}}}

\newcommand{\upchi}{\raise 2pt \hbox{$\chi$}}
\newcommand{\smallupchi}{\raise 1pt \hbox{\tiny $\chi$}}

\def\writefig#1 #2 #3 {\rlap{\kern #1 truecm \raise #2 truecm
\hbox{#3}}}

\newcommand{\caM}{{\mathcal M}}\newcommand{\caN}{{\mathcal N}}\newcommand{\caS}{{\mathcal S}}

\newcommand{\bbN}{{\mathbb N}}\newcommand{\bbR}{{\mathbb R}}\newcommand{\bbZ}{{\mathbb Z}}

\newcommand{\bsn}{{\boldsymbol n}}\newcommand{\bsr}{{\boldsymbol r}}\newcommand{\bsx}{{\boldsymbol x}}

\newcommand{\bsell}{{\boldsymbol\ell}}

\newcommand{\bspi}{{\boldsymbol\pi}}

\begin{document}

{\hfill\small Reference: Electr. J. Probab. 16, 1173-1192 (2011)}

{\hfill\small This version contains an extra appendix.} \vspace{3mm}

\title[Spatial random permutations and Poisson-Dirichlet law]{Spatial random permutations and Poisson-Dirichlet law of cycle lengths}

\author{Volker Betz and Daniel Ueltschi}

\thanks{Department of Mathematics, University of Warwick, Coventry, CV4 7AL, England. \\ v.m.betz@warwick.ac.uk, daniel@ueltschi.org}

\thanks{V.B.\ is supported by EPSRC grant EP/D07181X/1 and D.U. is supported in part by EPSRC grant EP/G056390/1.}

\maketitle

\begin{quote}
{\small
{\bf Abstract.}
We study spatial permutations with cycle weights that are bounded or slowly diverging. We show that a phase transition occurs at an explicit critical density. The long cycles are macroscopic and their cycle lengths satisfy a Poisson-Dirichlet law.}  

\vspace{1mm}

\vspace{1mm}
\noindent
{\footnotesize {\it 2010 Math.\ Subj.\ Class.:} 60K35, 82B26.}

\vspace{1mm}
\noindent
{\footnotesize Submitted 14 July 2010; accepted 11 May 2011.}

\end{quote}

\section{Introduction}

The structure for spatial permutations consists of a large box $\Lambda \subset \bbR^d$, a large number $N$ of points 
in $\Lambda$, and permutations of those points such that all permutation jumps remain small. The relevant parameter is 
the density $\rho = N/|\Lambda|$. In many models there is a critical density $\rho_{\rm{c}}$ that corresponds to a 
transition from a phase with only finite cycles (when $\rho \leq \rho_{\rm{c}}$) to a phase where a nonzero fraction of 
points belong to infinite cycles (when $\rho > \rho_{\rm{c}}$). The goal of the present article is twofold. First, we 
prove that such a transition occurs in a class of models of spatial random permutations with cycle weights. Second, we 
show that the cycle structure of infinite cycles satisfies a Poisson-Dirichlet law.

The main motivation for our models comes from the interacting Bose gas of quantum statistical mechanics. The possible  relevance of long permutation cycles to Bose-Einstein condensation was pointed out by Matsubara \cite{Mat} and Feynman \cite{Fey}. S\"ut\H o made important clarifications for the ideal Bose gas, showing in particular that long cycles are macroscopic \cite{Suto1,Suto2}. It is a notoriously difficult problem to prove Bose-Einstein condensation. Another problem, that is related but not subordinated, is to understand how the critical temperature is modified by particle interactions. In the recent article \cite{BU3}, we derived (non-rigorously) a model of spatial permutations where the original interactions between quantum particles have been replaced by cycle weights. The simplified model retains some features of the original model, as they have the same free energy to lowest order in the scattering length of the interaction potential. We then used the formula \eqref{crit dens} below for the critical density. The validity of this formula for the model with cycle weights is proved in the present article.

Models of spatial permutations are also attractive \textit{per se}. They have both specific and general features. One 
general feature that is especially striking is the Poisson-Dirichlet law for the distribution of cycle lengths. The 
literature on the subject is huge, see e.g.\ \cite{ABT,JKB,Han} for a sample. The Poisson-Dirichlet distribution is 
expected to make an appearance in other models with spatial structure and permutations such as the random stirring 
model \cite{Har,Toth}. This was proved recently by Schramm on the complete graph \cite{Sch}; see also Berestycki \cite
{Ber} for several useful observations and clarifications.

The models considered here are ``annealed" in the sense that spatial positions vary and they are integrated upon. 
Annealed models are both simpler and more relevant for the Bose gas. But the ``quenched" models, where the positions 
are fixed and chosen according to a suitable point process, look very interesting in probability theory. One 
conjectures that long cycles satisfy the same Poisson-Dirichlet law as in the annealed model --- the only difference 
being the critical density. This is supported by numerical evidence \cite{GRU,Kerl}. An unrelated but very interesting 
problem is the complete description of Gibbs states, involving crossing fluxes that depend on the boundary conditions. 
Such a description has been recently achieved by Biskup and Richthammer in the one-dimensional model \cite{BR}.

\section{Setting \& results}

The state space of the (annealed) model of spatial permutations with cycle weights is $\Omega_{\Lambda,N} = \Lambda^N 
\times \caS_N$, where $\Lambda \subset \bbR^d$ is a cubic box of size $L$, and $\caS_N$ is the symmetric group of 
permutations of $N$ elements. We denote by $|\Lambda|=L^d$ the volume of $\Lambda$. We equip $\Omega_{\Lambda,N}$ with 
the product of the Borel $\sigma$-algebra on $\Lambda^N$ and the discrete $\sigma$-algebra on $\caS_N$. We introduce a 
``Hamiltonian" and its corresponding Gibbs state. Namely, the Hamiltonian is a function $H : \Omega_{\Lambda,N} \to \bbR \cup \{\infty\}$ that we suppose of the form
\be
\label{Hamiltonian}
H(\bsx, \pi) = \sum_{i=1}^N \xi(x_i - x_{\pi(i)}) + \sum_{\ell\geq1} \alpha_\ell r_\ell(\pi).
\ee
Here, $\bsx = (x_1,\dots,x_N) \in \Lambda^N$ and $\pi \in \caS_N$. We always suppose that $\e{-\xi}$ is continuous with 
positive Fourier transform, and that it is normalized, $\int_{\bbR^d} \e{-\xi(x)} 
\dd x = 1$. Notice that $\xi$ is allowed to take the value $+\infty$, and that positivity of the Fourier transform implies that $\xi(x) = \xi(-x)$. The cycle weights $\alpha_1,\alpha_2,\dots$ are 
fixed parameters. Finally, $r_\ell(\pi)$ denotes the number of $\ell$-cycles in the permutation $\pi$.

Boundary conditions are not expected to play a prominent r\^ole here, and we therefore choose those that make proofs 
simpler. These are the ``periodized" boundary conditions, where we replace $\xi$ by the function $\xi_\Lambda$, defined 
by
\be
\e{-\xi_\Lambda(x)} = \sum_{z \in \bbZ^d} \e{-\xi(x - Lz)}.
\ee
The normalization assumption $\int_\Lambda \e{-\xi_\Lambda} = \int_{\bbR^d} \e{-\xi} = 1$ implies that $\e{-\xi_\Lambda(x)}$ is finite for at least almost every $x$.
When $\e{-\xi}$ has bounded support with diameter smaller than $L$ we recover the usual periodic boundary conditions.
We let $H_\Lambda$ be as $H$ in \eqref{Hamiltonian}, but with $\xi_\Lambda$ instead of $\xi$.
The Gibbs state is given by the probability measure
\be
\label{main measure}
{\rm Prob} (\dd \bsx, \pi) = \frac1{N! Y} \e{-H_\Lambda(\bsx,\pi)} \dd\bsx
\ee
on $\Omega_{\Lambda,N}$, where $\dd\bsx$ is the Lebesgue measure on $\Lambda^N$ and $Y$ is a suitable normalization, 
namely
\be
\label{normalization}
Y = \frac1{N!} \sum_{\pi\in\caS_n} \int_{\Lambda^N} \e{-H_\Lambda(\bsx,\pi)} \dd\bsx.
\ee

In typical realizations of the system, points are spread all over the space because of the Lebesgue measure that prevents 
accumulations. The lengths of permutation jumps $\| x_i - x_{\pi(i)}\|$ stay bounded uniformly in $\Lambda$ because of the 
jump weights $\e{-\xi(x_i - x_{\pi(i)})}$.  
The lengths of permutation cycles depend on the density of the system. For small density, points are far apart and jumps are 
unlikely, which typically results in small cycles. But as the density increases, points have more and more possibilities to 
hop, and a phase transition takes place where ``infinite" cycles appear. The cycle weights modify the critical density and 
also the distribution of cycle lengths, see below. The model is illustrated in Fig.\ \ref{fig spatial permut}.

\bfig
\centerline{\includegraphics[width=80mm]{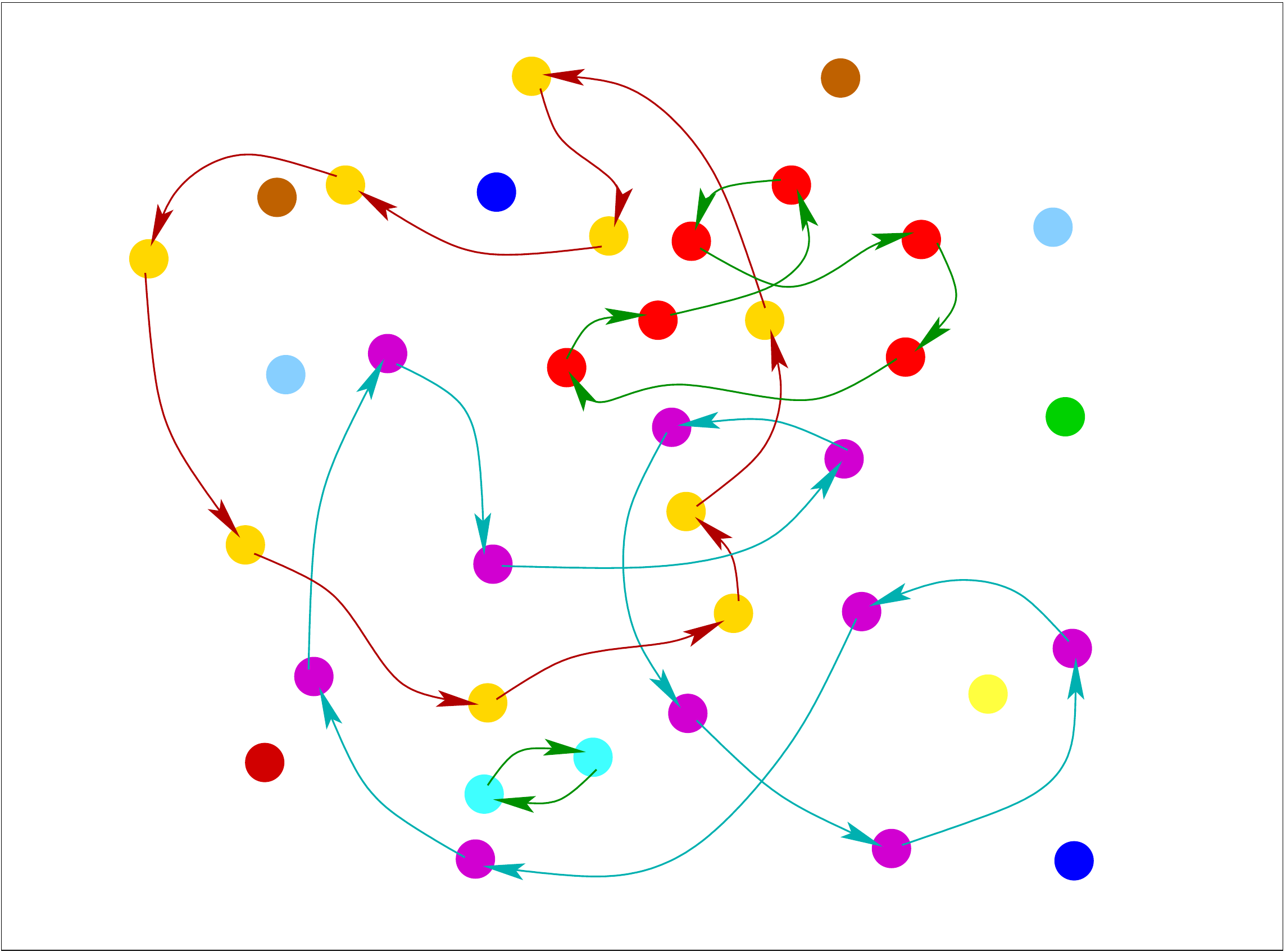}}
\caption{A typical realization of a spatial permutation. As $|\Lambda|,N \to \infty$, the jumps remain finite but the cycle lengths may diverge.}
\label{fig spatial permut}
\efig

This model arises naturally from the Feynman-Kac representation of the dilute Bose gas. The jump function is then $\xi(x) = \frac1{4\beta} \|x\|^2$ (plus a normalization constant), with $\beta$ the inverse temperature of the system. Notice that if the original quantum system has periodic boundary conditions, we get the periodized Gaussian function. Cycle weights were introduced in \cite{BU1} as a crude way to account for the particle interactions. But the calculations of \cite{BU3} suggest that the cycle weights can be chosen so that the model describes the Bose gas exactly in the dilute regime. We do not write here the precise formula for the weights, but we observe that they satisfy the asymptotic $\alpha_j = -\alpha (1 - O(j^{-1/5}))$, so that $\alpha_j$ converges as $j\to\infty$ fast enough for our purpose.

We are solely interested in properties of permutations and we introduce random variables that are functions on $\caS_N$ rather than $\Omega_{\Lambda,N}$. Let $\ell^{(1)}(\pi), \ell^{(2)}(\pi), \dots$ denote the cycle lengths in non-increasing order, repeated with multiplicities. We will prove that, above the critical density, the cycle lengths scale like $N$ and they converge in distribution to Poisson-Dirichlet. The latter is conveniently defined using the \textit{Griffiths-Engen-McCloskey distribution} GEM$(\theta)$, which is the distribution for
\[
\bigl( \, X_1, \, (1-X_1) X_2 \, , (1-X_1) (1-X_2) X_3 \, , \, \dots \bigr),
\]
where $X_1,X_2,\dots$ are i.i.d.\ beta random variables with parameter $(1,\theta)$; that is, ${\rm Prob}(X_i>s) = (1-s)^\theta$ for $0\leq s\leq 1$. The \textit{Poisson-Dirichlet distribution} PD$(\theta)$ is the law obtained by rearranging those numbers in non-increasing order. See \cite{ABT,JKB} for more information and background. In the sequel, we say that a sequence of random variables $Y^{(1)}_n,Y^{(2)}_n, \dots$ converges in distribution to Poisson-Dirichlet as $n\to\infty$ if, for each fixed $k$, the joint distribution of $Y^{(1)}_n,\dots,Y^{(k)}_n$ converges weakly to the joint distribution of the first $k$ random variables in Poisson-Dirichlet. This is denoted
\be
(Y^{(1)}_n,Y^{(2)}_n,\dots) \Rightarrow {\rm PD}(\theta).
\ee

As already mentioned, we make the important assumption that the jump function has nonnegative Fourier transform. This allows to define the ``dispersion relation" $\varepsilon(k)$, $k \in \bbR^d$, by the equation
\be
\label{eps k}
\e{-\varepsilon(k)} = \int_{\bbR^d} \e{-2\pi\ii k \cdot x} \e{-\xi(x)} \dd x.
\ee
Notice that $\eps(k)$ is real, $\eps(0) = 0$, and $\eps(k) > 0$ for all $k \neq 0$, and $\lim_{\|k\| \to \infty} \eps(k) = \infty$ (by Riemann-Lebesgue). In order to avoid pathological cases we assume that $\e{-\eps(k)}$ is uniformly continuous on $\bbR^d$. We also suppose that 
$\eps(k) \geq a \|k\|^\eta$ for small $k$, for some $a>0$ and $\eta<d$. 
It is easy to see that $\varepsilon(k)$ is always 
greater than $a\|k\|^2$ for small $k$, so the latter assumption always holds in dimensions $d>2$. Among possible jump functions 
other than 
Gaussians, let us mention $\e{-\xi(x)} = {\rm const} \, (|x|+1)^{-\gamma}$ with $1<\gamma<2$ in $d=1$, for which $\eta 
= \gamma-1$.
As for the cycle weights, we consider three cases:
\begin{itemize}
\item[(i)] $\lim_{j \to \infty} \alpha_j = \alpha$ with $\alpha > 0$, and $\sum_j |\alpha_j-\alpha| < \infty$.
\item[(ii)] $\lim_{j \to \infty} \alpha_j = \alpha$ with $\alpha \leq 0$, and $\sum_j \frac{1}{j}|\alpha_j-\alpha| < \infty$.
\item[(iii)] $\alpha_j = \gamma \log j$ with $\gamma > 0$.
\end{itemize}

We now introduce the fraction $\nu$ of points in infinite cycles. It is obvious that finite systems can only host finite cycles, so the definition of $\nu$ must involve the thermodynamic limit. Given a finite number $K$, let $\nu_K$ denote the fraction of points in cycles of length larger than $K$. Precisely,
\be
\label{def nu_K}
\nu_K = \liminf_{|\Lambda|,N \to \infty} E \Bigl( \frac1N \sum_{i : \ell^{(i)} > K} \ell^{(i)} \Bigr).
\ee
Here and in the sequel, the limit $|\Lambda|,N\to\infty$ means that both go to infinity while keeping the density $\rho = N/|\Lambda|$ fixed. This is the standard thermodynamic limit of statistical mechanics.
We then define
\be
\nu = \lim_{K\to\infty} \nu_K.
\ee
This limit exists since $(\nu_K)$ is decreasing and bounded. Let $\bar\nu_K$ denote the $\limsup$ of \eqref{def nu_K}. We expect that $\bar\nu_K = \nu_K$ but we do not prove it. On the other hand, we will prove in Section \ref{sec cycle lengths} that $\bar\nu_K$ also converges to $\nu$ as $K\to\infty$.

Next we introduce the \textit{critical density} by
\be
\label{crit dens}
\rho_{\rm{c}} = \sum_{j\geq1} \e{-\alpha_j} \int_{\bbR^d} \e{-j\varepsilon(k)} \dd k.
\ee
It follows from our assumptions that the critical density is finite. Indeed, the numbers $\e{-\alpha_j}$ are bounded, so $\rho_{\rm{c}}$ is bounded by the integral of a geometric series, $\int \frac1{\e{\eps(k)}-1}$, which is finite.

We propose now two theorems that confirm that $\rho_{\rm{c}}$ is indeed the critical density of the model, at least in several interesting situations. The formula \eqref{crit dens} is presumably valid beyond the cases treated in this article, but the precise extent of its validity is not clear.
The first theorem states that macroscopic cycles occur precisely above the critical density, and that they obey the Poisson-Dirichlet law.

\begin{theorem}
\label{thm Ewens}
Assume that $\alpha_j \to \alpha$ as described above. Then
\begin{itemize}
\item[(a)] the fraction of points in infinite cycles is given by
\[
\nu = \max\Bigl( 0, 1 - \frac{\rho_{\rm{c}}}\rho \Bigr);
\]
\item[(b)] when $\nu>0$, i.e.\ when $\rho>\rho_{\rm{c}}$, the cycle structure converges in distribution to Poisson-Dirichlet: As $|\Lambda|, N \to \infty$ we have 
\[
\Bigl( \frac{\ell^{(1)}}{\nu N}, \frac{\ell^{(2)}}{\nu N}, \dots \Bigr) \Rightarrow {\rm PD}(\e{-\alpha}).
\]
\end{itemize}
\end{theorem}

Such a law was already observed in absence of spatial structure, and when the cycle weights are constant. This case is known as the Ewens distribution, see e.g.\ \cite{Ewe,Han,ABT}. Results about weights that are asymptotically Ewens can be found in \cite{Lugo,BUV}. Spatial permutations with small cycle weights, i.e.\ when the limit is $\alpha=0$, were studied in \cite{BU2}.

The second theorem concerns cycle weights that diverge logarithmically  --- it is somehow the limit $\alpha\to\infty$ of Theorem \ref{thm Ewens}. Cycle weights have a striking effect as a single giant cycle occurs above the critical density! This is in accordance with a similar observation for non-spatial permutations \cite{BUV}.

\begin{theorem}
\label{thm giant}
Assume that $\alpha_j = \gamma \log j$ with $\gamma>0$. Then
\begin{itemize}
\item[(a)] the fraction of points in infinite cycles is given by
\[
\nu = \max\Bigl( 0, 1 - \frac{\rho_{\rm{c}}}\rho \Bigr);
\]
\item[(b)] when $\nu>0$, i.e.\ when $\rho>\rho_{\rm{c}}$, there is a single giant cycle that contains almost all points in infinite cycles: As $|\Lambda|,N \to \infty$ we have 
\[
\frac{\ell^{(1)}}{\nu N} \Rightarrow 1.
\]
\end{itemize}
\end{theorem}

The rest of this article is devoted to the proof of the results above. We reformulate the problem in the Fourier space in Section \ref{sec Fourier}, following S\"ut\H o \cite{Suto2}. The model involves a measure on occupation numbers of Fourier modes, and of random permutations of those numbers. In Section \ref{sec occupation numbers} we obtain information about occupation numbers using techniques of Buffet and Pul\'e \cite{BP}, and using certain estimates of our recent joint work with Velenik \cite{BUV}. Random permutations within each mode involve the cycle weights and are thus similar to those studied in \cite{BUV}. Combining all those results allows us to prove Theorems \ref{thm Ewens} and \ref{thm giant} in Section \ref{sec cycle lengths}.

\section{Random permutations and Fourier modes}
\label{sec Fourier}

The goal of this section is to introduce an alternative model of random permutations that involves Fourier modes, and that has the same marginal distribution on cycle lengths. Let $\Lambda^* = \frac1L \bbZ^d$ be the space dual to $\Lambda$ in the sense of Fourier theory.

\subsection{The marginal distribution of cycle lengths}

Recall that the cycle structure of a permutation $\pi \in \caS_N$ is the sequence of cycle lengths $\bsell = (\ell^{(1)}, \ell^{(2)}, \dots, \ell^{(m)})$, with $\ell^{(i)} \geq \ell^{(i+1)}$ and $\ell^{(m)} \geq 1$; the number of cycles $m$ depends on $\pi$, $1 \leq m \leq N$. Those numbers form a partition of $\{1,\dots,N\}$. Another way to write $\bsell$ is to introduce the occupation numbers $\bsr = (r_1,\dots,r_N)$, where $r_j = \#\{ i : \ell^{(i)} = j \}$. We always have
\be
\sum_{i=1}^m \ell^{(i)} = \sum_{j=1}^N j r_j = N.
\ee
One should not confuse the occupation numbers $\bsr$ with the occupation numbers $\bsn = (n_k)$ to be introduced later; they are not related in any direct way.

\begin{proposition}
\label{prop marginal}
The marginal of the probability measure \eqref{main measure} on occupation numbers is
\[
{\rm Prob}(\bsr) = \frac1Y \prod_{j=1}^N \frac1{r_j!} \Bigl( \frac{\e{-\alpha_j}}j \sum_{k \in \Lambda^*} \e{-j \eps(k)} \Bigr)^{r_j},
\]
with $Y$ the normalization of \eqref{normalization}.
\end{proposition}

\begin{proof}
The marginal probability on permutations is
\be
\begin{split}
{\rm Prob}(\pi) &= \frac1{N! Y} \int_{\Lambda^N} \e{-H_\Lambda(\bsx,\pi)} \dd\bsx \\
&= \frac1{N! Y} \int_{\Lambda^N} \e{-\sum_{i=1}^N \xi_\Lambda(x_i - x_{\pi(i)}) - \sum_{j\geq1} \alpha_j r_j(\pi)} \dd x_1 \dots \dd x_N.
\end{split}
\ee
We observe that integrals factorize according to permutation cycles. The contribution of a cycle of length $j$ is (with 
$y_{j+1} \equiv y_1$)
\be \label{convolution equality}
{}\e{-\alpha_j} \int_{\Lambda^j} \e{-\sum_{i=1}^j \xi_\Lambda(y_i - y_{i+1})} \dd y_1 \dots \dd y_j = \e{-\alpha_j} |
\Lambda| \sum_{z \in \bbZ^d} \bigl( \e{-\xi} \bigr)^{*j}(Lz).
\ee
To see the equality in \eqref{convolution equality}, we start with the right hand side. Using the definition of the convolution, 
writing $|\Lambda| = \int_{\Lambda} \dd y_1$, and shifting all the variables in the convolution integrals by $y_1$ 
gives 
\be
\begin{split}
&|\Lambda| \sum_{z \in \bbZ^d} \bigl( \e{-\xi} \bigr)^{*j}(Lz)  \\
&= \sum_{z \in \bbZ^d} \int_\Lambda \dd y_1
\int_{\bbR^{d(j-1)}} \dd y_2 \dots \dd y_j \e{-\xi(L z - y_2 + y_1)} \e{-\xi(y_2-y_3)} \dots \e{-\xi(y_j-y_1)}\\
&= \sum_{z_1 \in \bbZ^d} \int_\Lambda \dd y_1 \sum_{z_2, \ldots, z_j \in \bbZ^d} 
\int_{\Lambda^{j-1}} \dd y_2 \dots \dd y_j \e{-\xi(y_1-y_2 + L(z_1-z_2))} \e{-\xi(y_2-y_3 + L(z_2-z_3))} \\
& \hspace{4cm} \dots \e{-\xi(y_{j-1}-y_j + L(z_{j-1}-z_j))} - \e{-\xi(y_j-y_1 + L z_j)}.
\end{split}
\ee
The last equality is obtained by decomposing the domain of integration $\bbR^d$ into cubes $\Lambda + L z$ with 
$z \in \bbZ^d$ and then changing variables in the integrals so that all the boxes become centered at $0$. We now
change to summation index: $\tilde z_j = z_j$, and $\tilde z_i = z_i - z_{i+1}$ for $i<j$. It is easy to see that 
this is indeed a bijection on $(\bbZ^d)^j$. Summing over $\tilde z_i$ instead of $z_i$ now gives the left hand side of
\eqref{convolution equality}.

The Fourier transform of $(\e{-\xi})^{*j}$ is $\e{-j \eps(k)}$. The Poisson summation formula states that
\be
\label{Poisson}
\sum_{z \in \bbZ^d} f(Lz) = \frac1{|\Lambda|} \sum_{k \in \Lambda^*} \widehat f(k),
\ee
where $\widehat f$ is the Fourier transform of $f$, whose precise definition can be found in Eq.\ \eqref{eps k}. 
We then get
\be
{\rm Prob}(\pi) = \frac1{N! Y} \prod_{j=1}^N \Bigl( \e{-\alpha_j} \sum_{k \in \Lambda^*} \e{-j \eps(k)} \Bigr)^{r_j(\pi)}.
\ee
All permutations of a given cycle structure have the same probability, and there are
\be
\label{number in equiv class}
\frac{N!}{\prod_j j^{r_j} r_j !}
\ee
elements in the cycle structure defined by $\bsr$. We get the claim by multiplying the above probability by this number.
\end{proof}

\subsection{Decomposition of permutations according to Fourier modes}

We denote by $\bsn = (n_k)$ the occupation numbers indexed by $k \in \Lambda^*$, and by $\caN_{\Lambda,N}$ the set of occupation numbers such that $\sum_{k \in \Lambda^*} n_k = N$. Next, we introduce permutations that are also indexed by Fourier modes, $\bspi = (\pi_k)$. Let $\caM_{\Lambda,N}$ be the set of pairs $(\bsn,\bspi)$ where $\bsn \in \caN_{\Lambda,N}$ and $\bspi = (\pi_k)$ with $\pi_k \in \caS_{n_k}$ for each $k \in \Lambda^*$. We introduce a probability measure on $\caN_{\Lambda,N}$:
\be
\label{prob n}
{\rm Prob}(\bsn) = \frac1Y \prod_{k\in\Lambda^*} \e{-\eps(k) n_k} h_{n_k}
\ee
with
\be
\label{def h_n}
h_n = \frac1{n!} \sum_{\pi \in \caS_n} \e{-\sum_{j\geq1} \alpha_j r_j(\pi)},
\ee
and $h_0=1$.
We will check later that the normalization $Y$ is the same as given in \eqref{normalization}. Then we introduce the probability of a pair $(\bsn,\bspi)$ by
\be
\label{prob pair}
{\rm Prob}(\bsn,\bspi) = \frac1Y \prod_{k \in \Lambda^*} \frac1{n_k!} \e{-\eps(k) n_k - \sum_{j\geq1} \alpha_j r_j(\pi_k)}.
\ee
Notice that \eqref{prob n} is the marginal of \eqref{prob pair} with respect to $\bspi$. The conditional probability ${\rm Prob}(\bspi | \bsn)$, where $\pi_k \in \caS_{n_k}$ for all $k$, is given by
\be
{\rm Prob}(\bspi | \bsn) = \prod_{k\in\Lambda^*} \Bigl( \frac1{n_k! h_{n_k}} \e{-\sum_{j\geq1} \alpha_j r_j(\pi_k)} \Bigr).
\ee
That is, given $\bsn$, each $\pi_k$ is independent and distributed as nonspatial random permutations with cycle weights (see Eq.\ \eqref{prob nonspatial} below). Given $\bspi$, let $r_j = \sum_k r_j(\pi_k)$.

\begin{proposition}
The marginal of the probability measure \eqref{prob pair} with respect to $\bsr$ is identical to the marginal of the probability measure \eqref{main measure}.
\end{proposition}

\begin{proof}
We check that the marginal of \eqref{prob pair} gives the formula of Proposition \ref{prop marginal}. 
For this, let $\bsr$ be a collection of occupation numbers, and write $(r_{jk}):\bsr$ for the set of all 
integers $r_{jk}$ ($j \geq 1, k \in \Lambda^\ast$) such that $\sum_k r_{jk} = r_j$ for all $j$. Then,  
\be
\begin{split}
{\rm Prob}(\bsr) &= \frac1Y \sum_{(r_{jk}) : \bsr} \, \, \sumtwo{(\bsn,\bspi):}{r_j(\pi_k) = r_{jk}} 
\prod_{k\in\Lambda^*} \Bigl( \frac1{n_k!} \e{-\eps(k) n_k - \sum_j \alpha_j r_j(\pi_k)} \Bigr) \\
&= \frac1Y \sum_{(r_{jk}) : \bsr} \prod_{k\in\Lambda^*} \Bigl( \frac1{\prod_j j^{r_{jk}} r_{jk}!} \e{-\eps(k) \sum_j j r_
{jk} - \sum_j \alpha_j r_{jk}} \Bigr).
\end{split}
\ee
We have summed over $\pi_k$ that are compatible with $r_{jk}$, using the formula \eqref{number in equiv class} for the number of elements. The bracket above factorizes according to $j$. Using
\be
\prod_{k\in\Lambda^*} \frac{\e{-\alpha_j r_{jk}}}{j^{r_{jk}}} = \Bigl( \frac{\e{-\alpha_j}}j \Bigr)^{r_j},
\ee
and writing, for fixed $j \geq 1$, $(r_{jk}): r_j$ for the set of integers $r_{jk}$ with 
$\sum_k r_{jk} = r_j$, we get
\be
{\rm Prob}(\bsr) = \frac1Y \prod_{j\geq1} \biggl[ \Bigl( \frac{\e{-\alpha_j}}j \Bigr)^{r_j} \sum_{(r_{jk}) : r_j} \prod_{k\in\Lambda^*} \frac1{r_{jk}!} \e{-j \eps(k) r_{jk}} \biggr].
\ee
For each fixed $j$, the multinomial theorem gives
\be
\sum_{(r_{jk}) : r_j} \prod_{k\in\Lambda^*} \frac{\e{-j \eps(k) r_{jk}}}{r_{jk}!} = \frac1{r_j!} \Bigl( \sum_{k\in\Lambda^*} \e{-j \eps(k)} \Bigr)^{r_j}.
\ee
Then ${\rm Prob}(\bsr)$ is indeed given by the formula of Proposition \ref{prop marginal}. This also proves that $Y$ is the correct normalization that makes \eqref{prob n} and \eqref{prob pair} probability measures.
\end{proof}

\section{Properties of occupation numbers}
\label{sec occupation numbers}

We study in this section the probability measure of occupation numbers of Fourier modes, ${\rm Prob}(\bsn)$, that is defined in \eqref{prob n}. We show that the typical $\bsn$ has the following properties:
\begin{itemize}
\item $\frac{n_0}N = \max(0, 1 - \frac{\rho_{\rm{c}}}\rho)$;
\item $\frac1N \sum_{0<\|k\|<\delta} n_k$ is small when $\delta$ is small.
\item For all $\delta > 0$, $\frac{1}{N} \sum_{\|k\| \geq \delta} n_k 1_{n_k > M}$ is small when $M$ is large.
\end{itemize}

The behavior of the normalizations $h_n$ defined in \eqref{def h_n} play an important r\^ole. We assume that $h_n$ grows or vanishes at most polynomially, i.e., there are constants $C$ and $\kappa$ such that for all $n\geq1$,
\be
\label{estimates h_n}
(Cn)^{-\kappa} \leq h_n \leq (Cn)^\kappa.
\ee
We also need that certain ratios of $h_n$ be bounded. Precisely, for $s >1$, let
\be
\label{ratios h_n}
C(s) = \suptwo{m,n\geq1}{n/s < m < sn} \frac{h_m}{h_n}.
\ee
We assume that $C(s)$ is finite for any $s$.

Those properties have been verified in \cite{BUV} when $\alpha_j \to \alpha$ and $\alpha_j = \gamma \log j$. Indeed, one finds $h_n \sim n^{-r}$, with $r = \e{-\alpha}-1$ in the first case and $r = -1-\gamma$ in the second case. The results of the present article actually extend to other cycle weights,
as long as Eqs \eqref{estimates h_n} and \eqref{ratios h_n} hold true.

The radius of convergence of the generating function of $h_n$ is equal to 1. We have the following identity for all $\gamma>0$:
\be
\label{generating function}
\sum_{n\geq0} \e{-\gamma n} h_n = \exp \sum_{j\geq1} \tfrac1j \e{-\gamma j -\alpha_j}.
\ee
See \cite{BU2,Lugo} for the proof.

\subsection{Macroscopic occupation of the zero mode}

We use a strategy that is inspired by Buffet and Pul\'e in their study of the ideal Bose gas \cite{BP}. It consists in 
looking at the Laplace transform of the distribution of $\frac{n_0}N$. Let $Y(N)$ be the normalization of 
Eq.\ \eqref{normalization}. We now put the explicit dependence on $N$ because it is going to vary. Notice that $Y(N)$ 
also depends on $\Lambda$, but the domain is fixed throughout.

We have
\be
{\rm Prob}(n_0=j) = \frac{h_j}{Y(N)} \sumtwo{\bsn \in \caN_{\Lambda,N}}{n_0=j} \prod_{k\neq0} \e{-\eps(k) n_k} h_{n_k} = h_j \frac{\check Y(N-j)}{Y(N)},
\ee
with
\be
\check Y(N) = \sumtwo{\bsn \in \caN_{\Lambda,N}}{n_0=0} \prod_{k \in \Lambda^*} \e{-\eps(k) n_k} h_{n_k}.
\ee
Notice the relation
\be
Y(N) = \sum_{j=0}^N h_j \check Y(N-j).
\ee

We are often going to interchange infinite sums and products. Let $\caN$ be the set of finite sequences of integers. Notice that $\caN$ is countably infinite. The following lemma is easy to prove and sufficient for our purpose.

\begin{lemma}
\label{sum product trick}
Let $a(k,n)$ be a nonnegative function such that $a(k,0)=1$ for all $k$. Then
\[
\sum_{\bsn \in \caN} \prod_{k \geq 1} a(k,n_k) = \prod_{k \geq 1} \Big( \sum_{n \geq 0} a(k,n) \Big).
\]
(It is possible that both sides are infinite.)
\end{lemma}

\begin{proof}
Let $\ell(\bsn)$ be the index of the largest nonzero integer in $\bsn$. For every $m\geq1$ we have
\be
\sum_{\bsn\in\caN : \ell(\bsn) \leq m} \prod_{k=1}^{\ell(\bsn)} a(k,n_k) = \sum_{n_1,\dots,n_m \geq 0} \prod_{k=1}^m a(k,n_k) = \prod_{k=1}^m \Bigl( \sum_{n\geq0} a(k,n) \Bigr).
\ee
The left hand side and the right hand side are clearly increasing in $m$, and we obtain the lemma by letting $m\to\infty$.
\end{proof}

Lemma \ref{sum product trick} also holds for complex $a(k,n)$ under the assumption $\sum_{k \geq 0} \sum_{n \geq 1} |a(k,n)|$ be finite. This can be proved using dominated convergence, but it is not needed here.

We introduce a Riemann approximation to the critical density \eqref{crit dens} which will be useful in Proposition \ref{prop mu Lambda} below.
\be
\label{crit dens Riemann}
\rho_{\rm{c}}^{(\Lambda)} = \sum_{j\geq1} \e{-\alpha_j} \frac1{|\Lambda|} \sum_{k\neq0} \e{-j \eps(k)}.
\ee

\begin{lemma}
\label{crit dens Riemann convergence}
$\lim_{|\Lambda| \to \infty} \rho_{\rm c}^{(\Lambda)} = \rho_{\rm c}$.
\end{lemma}

\begin{proof}
Since $\e{-\eps(k)}$ is uniformly continuous, the Riemann sum $|\Lambda|^{-1} \sum_{k \neq 0} \e{- j \eps(k)}$ converges to $\int \e{- j \eps(k)} \, \dd k$ for each $j$. We need to show that the limit $|\Lambda| \to \infty$ can be interchanged with the sum over $j$ and we use dominated convergence. We sum separately over $\|k\|\leq1$ and $\|k\|>1$. For $\|k\|>1$ we use $\eps(k)>c>0$, so that
\be
\label{premiere borne}
\frac{1}{|\Lambda|} \sum_{\|k\| > 1} \e{-j \eps(k)} \leq \e{- j c/2} \frac{1}{|\Lambda|} \sum_{k \in \Lambda^\ast} \e{-\eps(k)} 
\ee
for $j \geq 2$. The latter sum is easily seen to converge using Eq. \eqref{Poisson}, so the right hand side is bounded by $C \e{-jc/2}$, which is summable. For $\|k\|\leq1$, we use $\eps(k) > a\|k\|^\eta$ with $a>0$. Since $\e{-ja\|k\|^\eta}$ is decreasing, we can estimate it with integrals, namely
\be
\label{sum integral bound}
\begin{split}
\frac1{|\Lambda|} \sum_{k \neq 0, \|k\| \leq 1} \e{-j \eps(k)} &\leq \frac1{|\Lambda|} \sum_{k \neq 0} \e{-aj \|k\|^\eta} \\
&\leq 2^d \int_{\bbR^d} \e{-aj \|k\|^\eta} \dd k = \frac{2^d d \pi^{d/2} \Gamma(\frac d\eta)}{\Gamma(\frac d2+1) \eta (aj)^{d/\eta}}.
\end{split}
\ee
The only relevant term in the upper bound is $j^{-d/\eta}$, which makes the sum over $j$ summable in Eq. \eqref{crit dens Riemann}. The claim follows from the dominated convergence theorem.
\end{proof}

Let $\caN_\Lambda = \cup_{N\geq0} \caN_{\Lambda,N}$ be the set of occupation numbers on $\Lambda^*$ where $\sum_k n_k$ is an arbitrary finite number.
Using Lemma \ref{sum product trick} we find
\be
\label{calc part fct}
\begin{split}
\check Z &= \sum_{N\geq0} \check Y(N) = \sumtwo{\bsn \in \caN_\Lambda}{n_0=0} \prod_{k \in \Lambda^*} \e{-\eps(k) n_k} h_{n_k} \\
&= \prod_{k \neq 0} \Bigl( \sum_{n\geq0} \e{-\eps(k) n} h_n \Bigr) = \exp \Bigl( \sum_{j\geq1} \tfrac1j \e{-\alpha_j} \sum_{k\neq0} \e{-j\eps(k)} \Bigr).
\end{split}
\ee
The last identity follows from \eqref{generating function} since $\inf_{k \neq 0} \eps(k) = c(\Lambda) > 0$.  $\check Z$ is finite by Lemma \ref{crit dens Riemann convergence}. This allows to introduce the following probability measure on $[0,\infty)$: 
\be
\label{mu Lambda}
\mu_\Lambda = \frac1{\check Z} \sum_{N\geq0} \check Y(N) \delta_{N/|\Lambda|}.
\ee
The motivation for $\mu_\Lambda$ is that the distribution of the occupation of the zero mode can be expressed as
\be
\begin{split}
{\rm Prob}(\tfrac{n_0}N \geq a) &= \sum_{j = \lceil aN \rceil}^N h_j \frac{\check Y(N-j)}{Y(N)} = \sum_{j=0}^{\lfloor (1-a)N \rfloor} h_{N-j} \frac{\check Y(j)}{Y(N)} \\
&= \frac{\int_0^{(1-a)\rho} h(|\Lambda| (\rho-s)) \, \dd\mu_\Lambda(s)}{\int_0^\rho h(|\Lambda| (\rho-s)) \, \dd\mu_\Lambda(s)}.
\end{split}
\ee
Here $h(x)$ can be any function interpolating the values $h_n$ at $n \in \bbN$, e.g.\ linear interpolation. We use the notation purely for convenience and 
will never evaluate $h(x)$ at non-integer points.

We now have all the elements that allow to state and to prove the key properties leading to the macroscopic occupation of the zero Fourier mode.

\begin{proposition}\hfill
\label{prop mu Lambda}
\begin{itemize}
\item[(a)] $\mu_\Lambda \to \delta_{\rho_{\rm{c}}}$ weakly as $|\Lambda| \to \infty$.
\item[(b)] Let $\lambda(\Lambda) \in \bbR$ such that $|\lambda(\Lambda)| \leq |\Lambda|^{\frac{1-\eta/d}2}$, then
\[
E_{\mu_\Lambda} \bigl( \e{\lambda(\Lambda) (X -\rho_{\rm{c}}^{(\Lambda)})} \bigr) \to 1.
\]
\end{itemize}
\end{proposition}

The parameter $\eta$ in the claim (b) is the one that appears in the condition for $\eps(k)$, see the paragraph after Eq.\ \eqref{eps k}.
The relevant aspect of the claim (b) is that the expectation is bounded uniformly in the domain even though $\lambda(\Lambda)$ diverges. 
Markov's inequality then gives the following concentration property for $|\Lambda|$ large enough, which will be used later:
\be
\label{concentration property}
\begin{split}
{\rm Prob}_{\mu_\Lambda}(|X-\rho_{\rm{c}}|>\epsilon) &\leq \e{-\epsilon |\Lambda|^{\frac{1-\eta/d}2}} E_{\mu_\Lambda} \Bigl( \e{|\Lambda|^{\frac{1-\eta/d}2} (X - \rho_{\rm{c}}^{(\Lambda)})} + \e{-|\Lambda|^{\frac{1-\eta/d}2} (X - \rho_{\rm{c}}^{(\Lambda)})} \Bigr) \\
&\leq 3 \e{-\epsilon |\Lambda|^{\frac{1-\eta/d}2}}.
\end{split}
\ee

\begin{proof}[Proof of Proposition \ref{prop mu Lambda}]
(a) follows from (b), see \eqref{concentration property}. For (b) we note that 
$a(k,n) = \e{-(\eps(k) - \lambda(\Lambda)/|\Lambda|)n} h_n$ fulfils the assumptions of Lemma \ref{sum product trick}, 
and thus   
\be
\label{calc transform}
\begin{split}
E_{\mu_\Lambda}(\e{\lambda(\Lambda) X}) &= \int_0^\infty \e{\lambda(\Lambda) s} \dd\mu_\Lambda(s) 
= \frac1{\check Z} \sum_{N\geq0} \check Y(N) \e{\lambda(\Lambda) N/|\Lambda|} \\
&= \frac1{\check Z} \sumtwo{\bsn \in \caN_\Lambda}{n_0=0} \prod_{k\neq0} \e{-(\eps(k) - \lambda(\Lambda)/|\Lambda|) n_k} h_{n_k}
\end{split}
\ee
Since $\inf_{k \neq 0} \eps(k) - \lambda(\Lambda)/|\Lambda| > 0$, \eqref{generating function} applies, and together with 
\eqref{calc part fct} we obtain
\be
E_{\mu_\Lambda}(\e{\lambda(\Lambda) X}) =  \exp \Bigl( \sum_{j\geq1} \tfrac1j \e{-\alpha_j} \sum_{k\neq0} 
\e{-j \eps(k)} (\e{j \lambda(\Lambda)/|\Lambda|} - 1) \Bigr) 
\ee
By \eqref{crit dens Riemann} and rearranging, we get 
\be
\label{Laplace transform}
E_{\mu_\Lambda} \bigl( \e{\lambda(\Lambda) (X-\rho_{\rm{c}}^{(\Lambda)})} \bigr) = \exp \Bigl( \sum_{j\geq1} \tfrac1j \e{-\alpha_j} \sum_{k\neq0} \e{-j\eps(k)} \bigl( \e{j\lambda(\Lambda)/|\Lambda|} - 1 - j \tfrac{\lambda(\Lambda)}{|\Lambda|} \bigr) \Bigr).
\ee
We show that the exponent vanishes as $|\Lambda|\to\infty$ using dominated convergence. We use $|\e{x} - 1 - x| \leq \frac12 x^2 \e{|x|}$ (which is easy to check using Taylor series) with $x = j \lambda(\Lambda)/|\Lambda|$.
The exponent in the right hand side of \eqref{Laplace transform} is bounded, in absolute value, by 
\be
\label{utile ici}
\tfrac12 \sum_{j \geq 1} \e{-\alpha_j} \frac{1}{|\Lambda|} \sum_{k \neq 0} \e{-\frac12 j \eps(k)} \Bigl[ \tfrac{j \lambda^2(\Lambda)}{|\Lambda|} \e{-j (\frac12 \eps(k) - |\lambda(\Lambda)|/|\Lambda|)} \Bigr].
\ee
Since $\varepsilon(k) \geq a |\Lambda|^{-\eta/d}$ and $|\lambda(\Lambda)|^2 \leq |\Lambda|^{1-\eta/d}$, the bracket is bounded above uniformly in $j$ for all $|\Lambda|$ large enough. The sum over $k$ has been estimated in \eqref{premiere borne} and \eqref{sum integral bound}. Since $\sum_j \e{-\alpha_j} j^{-d/\eta} < \infty$, we can interchange the limit $|\Lambda|\to\infty$ and the sum over $j$ by dominated convergence. The bracket in \eqref{utile ici} tends to 0 as $|\Lambda| \to \infty$, for all fixed $j$. It follows that \eqref{utile ici} converges to zero and \eqref{Laplace transform} converges to one. 
\end{proof}

\begin{proposition}
\label{prop large densities}
Suppose that $\rho > \rho_{\rm{c}}$. Then, in the thermodynamic limit $|\Lambda|,N \to \infty$,
\[
{\rm Prob}(\tfrac{n_0}N > a) \to \begin{cases} 1 & \text{if } a < 1 - \frac{\rho_{\rm{c}}}\rho, \\ 0 & \text{if } a > 1 - \frac{\rho_{\rm{c}}}\rho. \end{cases}
\]
\end{proposition}

\begin{proof}
We show that, for every $\epsilon>0$,
\be
{\rm Prob} \bigl( \bigl| \tfrac{n_0}N - \tfrac{\rho-\rho_{\rm c}}\rho \bigr| > \epsilon \bigr) \to 0.
\ee
Using the expression \eqref{mu Lambda} which involves the measure $\mu_\Lambda$, we write the probability as
\be
\label{basic split}
{\rm Prob} \bigl( \bigl| \tfrac{n_0}N - \tfrac{\rho-\rho_{\rm c}}\rho \bigr| > \epsilon \bigr) = \frac{J_- + J_+}{J_- + J_0 + J_+},
\ee
with 
\be
\begin{split}
&J_- = \int_0^{\rho_{\rm c} - \epsilon\rho} \frac{h(|\Lambda| (\rho-s))}{h(|\Lambda| (\rho - \rho_{\rm c}))} \dd\mu_\Lambda(s), \\
&J_0 = \int_{\rho_{\rm c} - \epsilon\rho}^{\rho_{\rm c} + \epsilon\rho} \frac{h(|\Lambda| (\rho-s))}{h(|\Lambda| (\rho - \rho_{\rm c}))} \dd\mu_\Lambda(s), \\
&J_+ = \int_{\rho_{\rm c} + \epsilon\rho}^\rho \frac{h(|\Lambda| (\rho-s))}{h(|\Lambda| (\rho - \rho_{\rm c}))} \dd\mu_\Lambda(s).
\end{split}
\ee
By \eqref{ratios h_n}, the ratios of $h(\cdot)$ are bounded above and below for all $0 < s < \rho - c$ uniformly in $\Lambda$. This shows that $J_- \to 0$ and that $J_0$ is bounded away from 0 as $|\Lambda|\to\infty$. As for $J_+$, we need to be careful with the integration over $s$ close to $\rho$. By \eqref{estimates h_n}, the ratio in $J_{+}$ is bounded by
\be
\frac{h(|\Lambda| (\rho-s))}{h(|\Lambda| (\rho - \rho_{\rm c}))} \leq \epsilon^{\kappa} (C|\Lambda|\rho)^{2\kappa},
\ee
which is valid provided $|\Lambda|\epsilon \geq 1$. Then
\be
J_{+} \leq \epsilon^{\kappa} (C|\Lambda|\rho)^{2\kappa} E_{\mu_{\Lambda}} (1_{s > \rho_{\rm c} + \epsilon \rho}),
\ee
and it goes to zero by Eq.\ \eqref{concentration property}.
\end{proof}

\subsection{No macroscopic occupation below the critical density}

For $\rho < \rho_{\rm{c}}$, the support of the Dirac measure to which the measures $\mu_\Lambda$ converge 
lies outside of the interval $[0,\rho]$, and the above argument fails. In its place, we use 
the formula 
\be \label{n_0 probability}
{\rm Prob}(n_0 \geq j) = \frac{1}{Y(N)} \sumtwo{\bsn \in \caN_{\Lambda,N}}{n_0 \geq j} \,\,
\prod_{k \in \Lambda^\ast} \e{- n_k \eps(k)} h_{n_k}
= \frac{Y(N-j,j)}{Y(N)},
\ee
where
\be
Y(m,j) = \sum_{\bsn \in \caN_{\Lambda,m}} 
\frac{h_{n_0+j}}{h_{n_0}}\prod_{k \in \Lambda^\ast} \e{- n_k \eps(k)} h_{n_k}. 
\ee
We apply summation by parts 
$E (f(X)) = f(0) + \sum_{j=1}^N \big[ f(j) - f(j-1) \big] P(X \geq j)$ to the function $\e{-\lambda n_0/N}$,
and we find
\be \label{expon RV}
\begin{split}
E(\e{-\lambda n_0/N}) &= 1 + \frac{(1-\e{\lambda/N})}{Y(N)}
\sum_{j=1}^{N} \e{-\lambda j/N} Y(N-j,j)\\
& = 
1 + \frac{\e{-\lambda}(1-\e{\lambda/N})}{Y(N)}
\sum_{j=0}^{N-1} \e{\lambda j/N} Y(j,N-j).
\end{split}
\ee
For the last equality, we used a change of summation index $j \to N-j$.

\begin{proposition}
\label{small densities}
Suppose that $\rho \leq \rho_{\rm{c}}$. In the thermodynamic limit $|\Lambda|$, $N \to \infty$,  
\[
{\rm Prob} (\tfrac{n_0}{N} > \delta) \to 0.
\]
for all $\delta > 0$. 
\end{proposition}

\begin{proof}
By \eqref{expon RV}, 
\be
| E(\e{\frac{\lambda n_0}{N}}) - 1|  \leq 
\frac{{\rm const}}{N} \left( \sum_{j=0}^{(1-\epsilon) N} \frac{Y(j,N-j)}{Y(N)} + 
\sum_{j=(1-\epsilon) N}^{N-1} \frac{Y(j,N-j)}{Y(N)} \right).
\ee
From \eqref{n_0 probability} it is obvious that $\frac{Y(j,N-j)}{Y(N)} \leq 1$, as it is equal to 
${\rm Prob}(n_0 \geq N-j)$. 
Thus the second term above, along with the prefactor $1/N$, is bounded by ${\rm const} \epsilon$. 
For the first term, we use the inequality 
\be
\sup_r \frac{h_{r+j}}{h_r} \leq C_0 (1+j)^{2\kappa}
\ee
for some constant $C_0$; this follows from \eqref{estimates h_n} and \eqref{ratios h_n}. Then, 
$Y(N,j) \leq Y(N) C_0 (1+j)^{2\kappa}$, 
and 
\be
\sum_{j=0}^{\lfloor (1-\epsilon) N \rfloor} \frac{Y(j,N-j)}{Y(N)} \leq C_0 
\sum_{j=0}^{\lfloor (1-\epsilon) N \rfloor} \frac{Y(j)}{Y(N)} (N - j + 1)^{2\kappa}
\ee
Now 
\be
\frac{Y(j)}{Y(N)} = \e{-|\Lambda| (q_{\Lambda}(j/|\Lambda|) - q_{\Lambda}(\rho))},
\ee
where $q_{\Lambda}(\rho) = -\frac{1}{|\Lambda|} \log Y(|\Lambda| \rho)$ is the finite volume free energy 
associated with the partition function $Y$. 
It was shown in \cite{BU2}, under conditions on the coefficients $\alpha_j$ that are more general than the present ones, that 
$q_{\Lambda}$ converges uniformly on compact intervals to a convex function $q$, and that
$\rho \mapsto q(\rho)$ is strictly decreasing for $\rho < \rho_{\rm{c}}$. 
Thus for each $\epsilon > 0$ there is $b_{\epsilon} > 0$ such that 
$q_\Lambda(j/|\Lambda|) - q_\Lambda(\rho) > b_\epsilon$ for all $|\Lambda|$ large enough, and all $j \leq (1-\epsilon) N$. So 
\be
\sum_{j=0}^{\lfloor (1-\epsilon) N \rfloor} \frac{Y(j,N-j)}{Y(N)} \leq C_0 \e{- b_\epsilon |\Lambda|} \sum_{j=0}^{\lfloor (1-\epsilon) N \rfloor} 
(N - j + 1)^{2\kappa} \leq C_0
\e{- b_\epsilon N / \rho} N^{2\kappa+1},
\ee
which converges to zero as $N \to \infty$. Since $\epsilon$ was arbitrary, we have shown that 
$E(\e{\lambda n_0/N}) \to 1$ for all $\lambda > 0$, which implies the claim. 
\end{proof}

\subsection{Occupation of nonzero modes}

We now turn to the modes $k \neq 0$. Recall the bound $C(s)$ for the ration of $h_n$ in Eq.\ \eqref{ratios h_n}.

\begin{lemma} \label{basic estimate}
For $0<\sigma <1$, $k \in \Lambda^\ast$, and $j \geq 0$, we have  
\[
{\rm Prob}(n_k \geq j) \leq C(\tfrac1\sigma)^2 \e{- j (1-\sigma) \eps(k)}.
\]
\end{lemma}

\begin{proof}
Using \eqref{ratios h_n}, we have
\be
\begin{split}
{\rm Prob}&(n_k \geq j) = \frac{1}{Y(N)} \sumtwo{\bsn \in \caN_{\Lambda,N}}{n_k \geq j} 
\prod_{k' \in \Lambda^\ast} \e{-n_{k'} \eps(k')} h_{n_{k'}} \\
& = \frac{1}{Y(N)} \sumtwo{\bsn \in \caN_{\Lambda,N-(1-\sigma) j}}{n_k \geq \sigma j} 
\e{-j (1-\sigma) \eps(k)} \frac{h_{n_k+(1-\sigma)j}}{h_{n_k}}
\prod_{k' \in \Lambda^\ast} \e{-n_{k'} \eps(k')} h_{n_{k'}} \\ 
& \leq \frac{C(1/\sigma)}{Y(N)} \e{-j (1-\sigma) \eps(k)} 
\sumtwo{\bsn \in \caN_{\Lambda,N-(1-\sigma) j}}{n_k \geq \sigma j} 
\prod_{k' \in \Lambda^\ast} \e{-n_{k'} \eps(k')} h_{n_{k'}}.
\end{split}
\ee
We indeed estimated $\frac{h_{n_k+(1-\sigma)j}}{h_{n_k}} \leq C(\frac{1-\sigma}\sigma) \leq C(1/\sigma)$. Since $\eps(k) \geq \eps(0) = 0$, we get an upper bound by replacing the constraint $n_k \geq \sigma j$ by $n_0 \geq \sigma j$. Then
\be
\begin{split}
{\rm Prob}&(n_k \geq j) \leq \frac{C(1/\sigma)}{Y(N)} \e{-j (1-\sigma) \eps(k)} 
\sumtwo{\bsn \in \caN_{\Lambda,N-(1-\sigma) j}}{n_0 \geq \sigma j} 
\prod_{k' \in \Lambda^\ast} \e{-n_{k'} \eps(k')} h_{n_{k'}} \\
& = \frac{C(1/\sigma)}{Y(N)} \e{-j (1-\sigma) \eps(k)} 
\sumtwo{\bsn \in \caN_{\Lambda,N}}{n_0 \geq j} \frac{h_{n_0 - (1-\sigma)j}}{h_{n_0}}
\prod_{k' \in \Lambda^\ast} \e{-n_{k'} \eps(k')} h_{n_{k'}} \\
& \leq \frac{C(1/\sigma)^2}{Y(N)} \e{-j (1-\sigma) \eps(k)} 
\sumtwo{\bsn \in \caN_{\Lambda,N}}{n_0 \geq j} 
\prod_{k' \in \Lambda^\ast} \e{-n_{k'} \eps(k')} h_{n_{k'}} \\
&\leq C(1/\sigma)^2 
\e{-j (1-\sigma) \eps(k)}.
\end{split}
\ee
\end{proof}

We now define three sets of occupation numbers, each of which will be shown to have measure 
close to one. Let $\tilde\nu = \max (0,1 - \frac{\rho_{\rm{c}}}{\rho})$; we will prove in the next section that $\tilde\nu = \nu$, but we do not know this yet. The sets are 
\be
\begin{split}
&A_\epsilon = \bigl\{ \bsn \in \caN_{\Lambda,N} : \bigl| \tfrac{n_0}N - \tilde\nu \bigr| < \epsilon \bigr\} \\
&B_{\epsilon,\delta} = \bigl\{ \bsn \in \caN_{\Lambda,N} : \sum_{0<\|k\|<\delta} n_k < \epsilon N \bigr\} \\
&C_{\epsilon,\delta,M} = \Bigl\{ \bsn \in \caN_{\Lambda,N} : \sumtwo{k \in \Lambda^*, \|k\|\geq\delta}{n_k \geq M} n_k < \epsilon N \Bigr\}.
\end{split}
\ee

\begin{proposition}
\label{prop typical occupation numbers}
For any $\rho>0$, we have in the thermodynamic limit $N, |\Lambda| \to \infty$:
\begin{itemize}
\item[(a)] For any $\epsilon>0$, ${\rm Prob} (A_\epsilon) \to 1$.
\item[(b)] For any $\epsilon>0$, there exists $\delta>0$ such that 
$\liminf {\rm Prob} (B_{\epsilon,\delta}) > 1-\epsilon$.
\item[(c)] For any $\epsilon, \delta>0$ there exists $M>0$ such that 
$
\liminf {\rm Prob} (C_{\epsilon,\delta,M}) > 1 - \epsilon
$.
\end{itemize}
\end{proposition}

\begin{proof}
The claim (a) immediately follows from Propositions \ref{prop large densities} and \ref{small densities}. For (b), we use Lemma \ref{basic estimate} with $\sigma = 1/2$ to get 
\be
E(n_k) = \sum_{i\geq1} {\rm Prob}(n_k \geq i) \leq 
C(2)^2 \sum_{i\geq 1} \e{-\eps(k) i/2 } = \frac{C(2)^2}{\e{\eps(k)/2}-1}.
\ee 
For every $\delta > 0$ we get, by Markov's inequality, 
\be
{\rm Prob}(B_{\epsilon,\delta}^{\rm{c}}) \leq \frac{C(2)^2}{\epsilon N} \sum_{0<\|k\|<\delta} 
\frac1{\e{\eps(k)/2}-1} \; \stackrel{N\to\infty}{\longrightarrow} \; 
\frac{C(2)^2}{\epsilon \rho} \int_{\|k\|<\delta} \frac{\dd k}{\e{\eps(k)/2}-1}.
\ee
By the assumption $\eps(k) \geq a \|k\|^\eta$ with $\eta < d$,
the integral is finite, and thus $\delta$ can be chosen so small that 
$\limsup {\rm Prob} (B_{\epsilon,\delta}^{\rm{c}}) < \epsilon$.

For (c), we define $F(\bsn) = \sum_{k \in \Lambda^\ast, \|k\| \geq \delta} n_k 1_{\{ n_k \geq M\}}$. Note that $
{\rm Prob} (C_{\eps,\delta,M}^{\rm{c}}) = {\rm Prob}(F \geq \eps N)$.
Now 
\be
\begin{split}
E(F/N) &= \frac{1}{N} \sum_{k \in \Lambda^\ast, \|k\| \geq \delta} E( n_k 1_{ n_k \geq M })\\
 &= \frac{1}{N} \sum_{k \in \Lambda^\ast, \|k\| \geq \delta} \Bigl( M {\rm Prob}(n_k \geq M) + 
\sum_{j > M} {\rm Prob}(n_k \geq j) \Bigr),
\end{split}
\ee
where the last equality is summation by parts. By Lemma \ref{basic estimate}, 
\be
\sum_{j > M} {\rm Prob}(n_k \geq j) \leq C(2)^2 \sum_{j=M+1}^\infty \e{-j \eps(k)/2} = C(2)^2 \e{-M \eps(k)/2} 
\frac{1}{\e{\eps(k)/2}-1}.
\ee
Define $c(\delta) = \inf_{\|k\| \geq \delta} \eps(k)$. Note that $c(\delta)>0$ for all $\delta > 0$ follows from the properties of $\eps(k)$ stated below \eqref{eps k}. Then, 
\be
\label{triangle}
E(F/N) \leq C(2)^2 M \e{- M c(\delta)/4} \frac{1}{N} \sum_{k \in \Lambda^\ast, \|k\| \geq \delta} 
\e{-M \eps(k)/4} \Big( 1 + \frac{1}{M (\e{\eps(k)/2} - 1)} \Big).
\ee
The sum above, along with the factor $1/N$, converges to a Riemann integral which is finite thanks to our 
conditions on $\eps(k)$. Therefore $\limsup E(F/N) \leq C_1 M \e{-M c(\delta)/4}$, and Markov's inequality implies that
$\limsup {\rm Prob} (C^{\rm{c}}_{\epsilon,\delta,M}) \leq \frac{M}{\epsilon} \e{-M c(\delta)/4}$. Choosing $M$ large enough 
for given $\epsilon, \delta$ proves the claim. 
\end{proof}

\section{Cycle lengths of spatial permutations}
\label{sec cycle lengths}

We now prove the claims of Theorems \ref{thm Ewens} and \ref{thm giant}, starting with the fraction $\nu$ of points in infinite cycles. We denote by ${\rm Prob}_n(\pi)$ the probability of a permutation $\pi \in \caS_n$ in the nonspatial model with cycle weights. That is,
\be
\label{prob nonspatial}
{\rm Prob}_n(\pi) = \frac1{h_n n!} \prod_{j\geq1} \e{-\alpha_j r_j(\pi)}
\ee
with $h_n$ the normalization defined in \eqref{def h_n}. We also write $E_n$ for the corresponding expectation. We keep the notation ${\rm Prob}$, $E$ for probability and expectation with respect to the spatial model.

Recall that we defined $\tilde\nu = \max( 0, 1 - \frac{\rho_{\rm{c}}}\rho)$.

\begin{proposition}
\label{prop nu}
Under the assumptions of Theorems \ref{thm Ewens} or \ref{thm giant}, we have $\nu = \tilde\nu$. 
\end{proposition}

\begin{proof}
We use the Fourier modes decomposition of Section \ref{sec Fourier}. Recall that $\bspi = (\pi_k)$, $r_{jk} = r_j(\pi_k)$, and $r_j = \sum_k r_{jk}$. We have, for fixed $K>1$, 
\bm
\label{c'est utile}
E \Bigl( \frac1N \sum_{i: \ell^{(i)} > K} \ell^{(i)} \Bigr) = E \Bigl( \frac1N \sum_{j>K} j r_j \Bigr) \\
= E \Bigl( \frac1N \sum_{j>K} j r_{j0} \Bigr) + E \Bigl( \frac1N \sum_{0<\|k\|<\delta} \sum_{j>K} j r_{jk} \Bigr) + E \Bigl( \frac1N \sum_{\|k\| \geq \delta} \sum_{j>K} j r_{jk} \Bigr).
\end{multline}
The first term of the right-hand side is equal to
\be
E \Bigl( \frac1N \sum_{j>K} j r_{j0} \Bigr) = \sum_{n\geq0} \frac nN {\rm Prob}(n_0=n) E_{n} \Bigl( \frac1n \sum_{j>K} j r_j \Bigr).
\ee
It follows from Proposition \ref{prop typical occupation numbers} (a) that $\frac{n_0}N \to \tilde\nu$ in probability 
as $|\Lambda|,N \to \infty$. In addition, we have
\be
E_n \Bigl( \frac1n \sum_{j>K} j r_j \Bigr) = {\rm Prob}_n(\ell_1 > K),
\ee
where $\ell_1$ is the length of the cycle that contains the index 1. It was shown in \cite{Lugo,BUV} that the latter converges to 1 as $n\to\infty$. We have thus proved that, for any finite $K$,
\be
\lim_{|\Lambda|,N \to \infty} E \Bigl( \frac1N \sum_{j>K} j r_{j0} \Bigr) = \tilde\nu.
\ee
The second term in the right-hand side of \eqref{c'est utile} is less than $E(\frac1N \sum_{0<\|k\|<\delta} n_k)$ and this is as small as we want by choosing $\delta$ small, see Proposition \ref{prop typical occupation numbers} (b). The last term is less than
\[
E \Bigl( \frac1N \sum_{\|k\| \geq \delta} n_k 1_{n_k>K} \Bigr).
\]
For any $\delta>0$, this can be made small by choosing $K$ large, see Proposition \ref{prop typical occupation numbers} (c). This shows that both $\nu_K$ and $\bar\nu_K$ converge to $\tilde\nu$ as $K\to\infty$.
\end{proof}

The next step is to prove that the distribution of cycle lengths of nonspatial weighted random 
permutations is asymptotically equal to Poisson-Dirichlet.

\begin{proposition}
\label{prop PD}
Assume that $\alpha_j \to \alpha$ as in Theorem \ref{thm Ewens}. Then, under the probability measure \eqref{prob nonspatial},
\[
\Bigl( \frac{\ell^{(1)}}n, \dots, \frac{\ell^{(m)}}n \Bigr) \Rightarrow {\rm PD}(\e{-\alpha}).
\]
\end{proposition}

\begin{proof}
Let us order the cycles of a permutation $\pi$ according to some rule, such as their smallest element. That is, the first cycle is the one that contains the index 1; the second cycle is the one that contains the smallest element that is not already in the first cycle; and so on... Let $\ell_1,\ell_2,\dots$ be the cycle lengths with respect to this order. We prove that
\[
\Bigl( \frac{\ell_1}n, \frac{\ell_2}{n-\ell_1}, \dots, \frac{\ell_m}{n-\ell_1-\dots-\ell_{m-1}} \Bigr)
\]
converges (in distribution) to i.i.d.\ beta random variables with parameters $(1,\e{-\alpha})$. It then immediately follows that $(\frac{\ell_1}n,\dots,\frac{\ell_m}n)$ converges to ${\rm GEM}(\e{-\alpha})$, and that $(\frac{\ell^{(1)}}n,\dots,\frac{\ell^{(m)}}n)$ converges to ${\rm PD}(\e{-\alpha})$. We proceed by induction on $m$. The case $m=1$ is just the law for $\frac{\ell_1}n$, whose convergence to the beta random variable was proved in \cite{Lugo,BUV}. For $m\geq1$, let
\be
A = \Bigl\{ (c_1,\dots,c_m) \in \{1,\dots,n\}^m : \tfrac{c_1}n \leq a_1, \dots, \tfrac{c_m}{n - c_1 - \dots - c_{m-1}} \leq a_m \Bigr\}.
\ee
Then
\bm
\label{c'est trop long}
{\rm Prob}_n \Bigl( \bigl( \tfrac{\ell_1}n, \dots, \tfrac{\ell_{m+1}}{n - \ell_1 - \dots - \ell_m} \bigr) \in A \times [0,a_{m+1}] \Bigr) \\
= \sum_{(c_1,\dots,c_m) \in A} {\rm Prob}_n \Bigl( \ell_1 = c_1, \dots, \ell_m = c_m \Bigr) \\
{\rm Prob}_n \Bigl( \tfrac{\ell_{m+1}}{n - \ell_1 - \dots - \ell_m} \leq a_{m+1} \, \big| \, \ell_1 = c_1, \dots, \ell_m = c_m \Bigr).
\end{multline}
It is not hard to check that
\bm
\label{aussi long}
{\rm Prob}_n \Bigl( \tfrac{\ell_{m+1}}{n - \ell_1 - \dots - \ell_m} \leq a_{m+1} \, \big| \, \ell_1 = c_1, \dots, \ell_m = c_m \Bigr) \\
= {\rm Prob}_{n-c_1-\dots-c_m} \Bigl( \tfrac{\ell_1}{n-c_1-\dots-c_m} \leq a_{m+1} \Bigr).
\end{multline}
In addition, any element $(c_1,\dots,c_m)$ in $A$ satisfies
\be
n - c_1 - \dots - c_m \geq n \prod_{i=1}^m (1-a_i).
\ee
It follows that \eqref{aussi long} converges in probability to the beta measure of $[0,a_{m+1}]$ uniformly in $c_1,\dots,c_m$. Then \eqref{c'est trop long} converges to a product of beta measures of the set $\times_{i=1}^{m-1} [0,a_i]$ by the induction hypothesis.
\end{proof}

Finally, we relate the distribution of long cycles of the spatial model to that of nonspatial random permutations. Let
\be
\label{def A}
A = [a_1,b_1] \times \dots \times [a_m,b_m] \subset (0,1)^m.
\ee

\begin{proposition}
\label{prop equivalence}
If $\nu>0$, we have for any $m\geq1$,
\[
\lim_{|\Lambda|,N\to\infty} {\rm Prob} \Bigl( \Bigl( \frac{\ell^{(1)}}{\nu N}, \dots, \frac{\ell^{(m)}}{\nu N} \Bigr) \in A \Bigr) = \lim_{n\to\infty} {\rm Prob}_n \Bigl( \Bigl( \frac{\ell^{(1)}}{n}, \dots, \frac{\ell^{(m)}}{n} \Bigr) \in A \Bigr).
\]
\end{proposition}

\begin{proof}
Let us write $\ell^{(1)}_k$ for the length of the longest cycle of the permutation $\pi_k$ corresponding to $k \in \Lambda^\ast$. 
We clearly have
\be
{\rm Prob} \Bigl( \sup_{k\neq0} \frac{\ell^{(1)}_k}N > \epsilon \Bigr) \leq {\rm Prob} \Bigl( \sup_{k\neq0} \frac{n_k}N > \epsilon \Bigr).
\ee
It follows from Proposition \ref{prop typical occupation numbers} (b) and (c) that the right-hand side vanishes in the limit $|\Lambda|,N \to \infty$. The zero Fourier mode is consequently the only one that matters, i.e.\
\begin{align}
\lim_{|\Lambda|,N \to \infty} {\rm Prob} \Bigl( \Bigl( \frac{\ell^{(1)}}{\nu N}, \dots, & \frac{\ell^{(m)}}{\nu N} \Bigr) \in A \Bigr) = \lim_{|\Lambda|,N \to \infty} {\rm Prob} \Bigl( \Bigl( \frac{\ell^{(1)}_0}{\nu N}, \dots, \frac{\ell^{(m)}_0}{\nu N} \Bigr) \in A \Bigr) \nn\\
&= \lim_{|\Lambda|,N\to\infty} {\rm Prob} \Bigl( \Bigl( \frac{\ell^{(1)}_0}{n_0}, \dots, \frac{\ell^{(m)}_0}{n_0} \Bigr) \in A \Bigr).
\end{align}
The last identity follows from Proposition \ref{prop typical occupation numbers} (a). Since $n_0 \to \infty$ as $|\Lambda|,N\to\infty$, the last term converges to the asymptotic joint probability of the $m$ largest cycles in nonspatial random permutations with cycle weights.
\end{proof}

Theorem \ref{thm Ewens} clearly follows from Propositions \ref{prop nu}, \ref{prop PD}, and \ref{prop equivalence}. Theorem \ref{thm giant} follows from Propositions \ref{prop nu} and \ref{prop equivalence}, and from the fact that $\frac{\ell_1}n \Rightarrow 1$ for random permutations with cycle weights of the form $\e{-\alpha_j} = j^{-\gamma}$ with $\gamma>0$, see \cite{BUV}. Notice that Proposition \ref{prop equivalence} is trivial here for $m \geq 2$, as both sides of the identity converge to zero.

\medskip
\noindent
\textit{Acknowledgments:} It is a pleasure to thank Nathana\"el Berestycki, Nick Ercolani, Alan Hammond, James Martin, and Yvan Velenik for many enlightening discussions. We are also grateful to a referee for numerous suggestions that helped us to improve the clarity of the presentation.

\appendix

\section*{Addendum}

As Antal Jarai has pointed out to us, the proofs given in this article are not complete: We use the Poisson summation formula without assuming decay properties of the function; and our assumptions are not sufficient to ensure convergence of the Riemann approximations to integrals.

We actually need to make two additional assumptions on the function $\e{-\xi(x)}$. 
These assumptions seem to hold for all reasonable choices of $\xi$, and in fact we have not 
been able to find an example of a function $\e{-\xi}$ that fullfills our original assumptions (in particular to be positive with
positive Fourier transform), but not the new ones. The main statements of the article stand unchanged.

We first discuss the Poisson formula. Rather than assuming extra decay property, we use positivity of the function and its Fourier transform.\\[1mm]

{\bf First additional assumption.} Suppose that
\[
\e{-\xi_\Lambda(x)} = \sum_{z \in \mathbb Z^{d}} \e{-\xi(x-L z)}
\]
is a continuous function of $x$, for all $L$ large enough. It follows that $(\e{-\xi_{\Lambda}})^{*j}$ is also continuous and integrable (Young's inequality). Next, let us recall the following Fourier inversion formula, that is valid for any continuous function $f$ on $[-1,1]^d$:
\[
f(x) = \lim_{\eta \to 0} \sum_{k \in \mathbb Z^d} \hat f(k) \e{-\pi \eta |k|^2} \e{2 \pi \ii k \cdot x},
\]
for all $x$ in the torus, and where 
\[
\hat f(k) = \int_{[-1,1]^d} f(x) \e{-2 \pi \ii k\cdot x} \dd x.
\]
This gives, after rescaling,
\[
(\e{-\xi_\Lambda})^{*j}(x) = \frac{1}{|\Lambda|} \lim_{\eta \to 0} \sum_{k \in \mathbb Z^d/L} \e{- j \varepsilon(k)} \e{- \pi \eta |k|^2} \e{2 \pi \ii k \cdot x}.
\]
For $x=0$ we can use the monotone convergence theorem to take the limit $\eta \to 0$, and we find 
\[
(\e{-\xi_\Lambda})^{*j}(0) \equiv \sum_{z \in \mathbb Z^d} (\e{-\xi})^{*j}(L z) = \frac{1}{|\Lambda|} \sum_{k \in \Lambda^*} \e{-j \varepsilon(k)}.
\]
This allows to get Eq.\ (3.6), from Eq.\ (3.3).
\\[1mm]

Next, in the proof of Lemma 4.2, we mistakenly assumed that Riemann sums of uniformly continuous functions always converge to integrals. There are indeed counter-examples, though such functions may not be positive with a positive Fourier transform.
\\[1mm]

{\bf Second additional assumption.} Suppose that
\[
\lim_{L\to\infty} \frac{1}{|\Lambda|} \sum_{k \in \Lambda^*} \e{-\varepsilon(k)} = \int_{\mathbb R^{d}} \e{-\varepsilon(k)} \dd k.
\]
It follows from continuity of $\e{-\varepsilon(k)}$ that
\[
\lim_{L\to\infty} \frac{1}{|\Lambda|} \sum_{k \in \Lambda^* \cap [-R,R]^{d}} \e{-j \varepsilon(k)} = \int_{[-R,R]^{d}} \e{-j \varepsilon(k)} \dd k
\]
for all $j\geq1$, and monotone convergence allows to take $R\to\infty$.
\\[2mm]

Finally, our proof that (4.18) tends to 0 uses that $\e{-\frac12 \varepsilon(k)}$ is Riemann integrable, which we have not assumed. But it is not hard to repair this. First, the term $j=1$ in the exponent of (4.17) is easily seen to go to 0. The sum over $j\geq2$ can be bounded by (4.18) with $j\geq2$. The bracket is less than a constant. We get the bound
\[
C \sum_{j\geq1} \frac1{|\Lambda|} \sum_{k\neq0} \Bigl( \e{-j \varepsilon(k)} + \e{-(j+\frac12) \varepsilon(k)} \Bigr) \leq 2C \sum_{j\geq1} \frac1{|\Lambda|} \sum_{k\neq0} \e{-j \varepsilon(k)},
\]
which is finite, see the proof of Lemma 4.2. Notice that the constant $C$ above also depend on the weights $\e{-\alpha_{j}}$, but these are bounded.
\\[1mm]

We are very grateful to Antal Jarai for pointing out the gaps in the proofs.

\end{document}